    \documentclass[final,3p,times]{elsarticle}

 \usepackage{graphicx}

\usepackage{amsmath,amssymb}
\usepackage{color}

\newcommand{\be}{\begin{equation}}
\newcommand{\ee}{\end{equation}}

\begin{document}

\begin{frontmatter}



\title{A Finite Difference Method with Non-uniform Timesteps for Fractional Diffusion Equations}


\author{Santos B. Yuste \corref{cor1}}
\ead{santos@unex.es}
\ead[url]{http://www.unex.es/eweb/fisteor/santos/sby.html}
\author{Joaqu\'{\i}n Quintana-Murillo}

\cortext[cor1]{ Corresponding author}


\address{Departamento de F\'{\i}sica, Universidad de Extremadura,
E-06071 Badajoz, Spain}

\begin{abstract}
An implicit finite difference method with non-uniform timesteps for solving the fractional diffusion equation in the Caputo form is proposed.
The method allows one to build adaptive methods where the size of the timesteps is adjusted to the behaviour of the solution in order to keep  the numerical errors small without the penalty of a huge computational cost.
The method is unconditionally stable and convergent. In fact, it is shown  that consistency and stability implies convergence for a rather general class of fractional finite difference methods to which the present method belongs.
The huge computational advantage  of adaptive methods against fixed step methods for fractional diffusion equations is illustrated by solving the problem of the dispersion of a flux of subdiffusive particles stemming from a point source.

\end{abstract}


\end{frontmatter}


\section{Introduction}
\label{sec:Intro}

Fractional calculus, a very old field of Mathematics dating back to the time of Leibniz \cite{Oldham}, has recently become a research area of growing interest due  mainly  to its surprisingly broad range of applications  in Physics, Engineering, Chemistry, Biology, Economics \cite{HilferEd,MetzlerKlafterPhysRep,PhysToday,MetzlerKlafterJPA04,Kilbas:06,Anotrans:08,Magin:08}.
In particular, fractional calculus is a key tool in the study of some anomalous diffusion processes, which, for example,  are abundant in biological environments where the presence of traps and obstacles often leads to mean square displacements that grow sublinearly with time \cite{DixVerkman08,HenryLanglands08}. A well-known model for describing this kind of subdiffusion processes is the so-called Continuous Time Random Walk model in which waiting time distributions between successive steps have a power-law tail. A convenient property of this mesoscopic model is that it leads to macroscopic fractional diffusion equations that describe how the concentration of walkers evolves in space and time \cite{MetzlerKlafterPhysRep,BarMetKlaPRE00}.

Some analytical methods of solution of these equations are known (method of images, separation of variables, integral transform methods,\ldots) \cite{MetzlerKlafterPhysRep,Kilbas:06,YusAcePhysica04,Sayed,Jafari,Ray}. However, as is also the case for non-fractional problems, often it is not possible to find an (at least useful) analytic solution and one has to resort to numerical methods.  The proposal, development, and analysis of numerical methods to solve fractional differential equations is at present a  quite active field of research, and many and varied methods have been considered
\cite{GorenfloNLD02,
YusteSIAM05,YusJCP06,
LynchJCP03,
MeerschaertTadjeranJCP4,
Chen07,Sun08,
Podlubny09,Cui09,BrunnerJCP10,SkovranekFDA10,
Mustapha11}.
Among these methods, finite difference methods are particularly convenient, and  consequently have been extensively studied. However, in almost all the cases considered the timesteps are constant (a recent exception is the work by Skovranek et al. \cite{SkovranekFDA10}).   Certainly, in this way these methods are simpler. On the other hand, methods that allow  variable timesteps have the great advantage that they can be converted into adaptive methods in which the size of the timesteps is chosen according to the behaviour of the solution. For example, one can choose  small timesteps when the solution changes quickly in order to avoid the typical problem of methods with constant timesteps of passing over, or scarcely sampling, those regions where the solution has this behaviour. Besides being more reliable, adaptive methods are usually faster because they can employ large timesteps whenever the solution changes smoothly. In particular, one could dynamically adjust the size of the timestep so that the error is smaller than a prefixed value. The aim of this paper is to set the basis to extend adaptive methods to fractional diffusion equations by presenting a finite difference method with variable timesteps.

The use of adaptive finite difference methods to fractional equations is especially relevant because, for these equations, the number of operations required to calculate the numerical solution at the $n$th timestep scales as  $n^2$. In comparison,  this number simply scales as $n$ for normal equations. The $n^2$ behaviour is due to the requirement of using the values of the numerical solution for all the previous times at which the solution was calculated (because the fractional derivatives are non-local integro-differential operators), which means that the number of calculations is roughly proportional to $\sum_{i=1}^n i\sim n^2$. This makes these methods slow and hugely memory demanding when $n$ is large.  The problem is so acute that a great deal of effort  has been devoted to overcoming it. Two main approaches have been explored. The most obvious one consists of increasing the ability of the numerical method to evaluate the non-local operators with large timesteps without loss of accuracy \cite{DiethelmND02,DiethelmCMAME05}. A less straightforward approach is that of the nested meshes method which is based on some scaling properties of the fractional derivatives \cite{FordNA01}. The number of operations in this method scales only as $n \ln n$.
It should be noted that these two approaches are compatible with the method with non-uniform temporal meshes to be discussed in this paper.

The stability of a finite difference method is, from a practical point of view, the key property one has to assess because it is not useful numerically if one does not know under what circumstances the method is unstable. The stability of the present method is analyzed by means of the von Neumann (or Fourier) procedure particularized for fractional diffusion equations  \cite{YusteSIAM05,YusJCP06,Murio08}. It is found that the method is unconditionally stable regardless of the temporal non-uniform mesh employed.  Besides, it is proved in a rather general form that under relatively weak conditions  the consistency and stability of a fractional finite difference method implies its convergence, i.e., it implies that the solution of the continuous equation is recovered from the finite difference solution when the size of the discretization mesh goes to zero.

We shall consider the fractional diffusion equation in the Caputo form
\begin{subequations}
\label{FDEdef}
\begin{equation}
\partial u=F
 \label{FDEeq}
\end{equation}
with
\begin{equation}
\partial \equiv \frac{\partial^\gamma}{\partial t^\gamma}-
  K \frac{\partial^2}{\partial x^2} ,
 \label{FDEoperator}
\end{equation}
\end{subequations}
and where
\begin{equation}
\label{CaputoRL2b}
  \frac{\partial^\gamma}{\partial t^\gamma} y(t)\equiv
\frac{1}{\Gamma(1-\gamma)} \int_0^t d\tau  \frac{1}{(t-\tau)^{\gamma}}
\frac{dy(\tau)}{d\tau},\quad 0<\gamma<1
\end{equation}
is the Caputo fractional derivative and $F(x,t)$ is a given source term.
It should be noted that our procedure can be extended to other equations with other terms (fractional Fokker-Planck equations, fractional diffusion-wave equations where $1<\gamma<2$) or even with other fractional operators such as the Riemann-Liouville derivative.
 In some  diffusion problems $u(x,t)$ represents the probability density of finding a particle at $x$ at time $t$. When the particle starts at zero at time zero, the solution $u(x,t)$ of $\partial u=0$ is the propagator (or Green's function) with $u(x,0)$ given by  Dirac's delta function, $\delta(x)$, and  $\langle x^2 \rangle \sim 2 K t^\gamma/\Gamma(1+\gamma)$ is the corresponding mean square displacement of the particle for large $t$. For this reason the parameter $K$ is sometimes called the anomalous diffusion coefficient  and $\gamma$  the anomalous diffusion exponent. The operator $d^\gamma/dt^\gamma$ is the usual first order derivative when  $\gamma=1$

The fractional diffusion equation is also usually written in terms of the fractional Riemann-Liouville derivative \cite{MetzlerKlafterPhysRep}. From a practical perspective, the two representations are different ways of writing the same equation as they are equivalent under fairly general conditions\cite{YusteQuintanaPS09}.
However, the numerical methods obtained for each representation are formally different, which leads to different ease of use and efficiency of the numerical algorithm. For example, in \cite{YusteQuintanaPS09} it was proved that the method of Gorenflo et al. \cite{GorenfloNLD02} for the diffusion equation in the Caputo form and the method of  Yuste and Acedo in \cite{YusteSIAM05} for the diffusion equation in the Riemann-Liouville form were equivalent only if the Gr\"{u}nwald-Letnikov discretization (or BDF1 formula) is used in both methods.

The difference scheme used in this paper is obtained by discretizing (i) the Laplacian with the three-point centred formula and (ii) the Caputo derivative with a generalization of the L1 formula  to non-uniform meshes.  It is second-order accurate in the spatial mesh size and first-order accurate in the timestep size. For a uniform temporal mesh, it becomes the numerical scheme discussed by Liu et al. \cite{LiuANZIAM06} and Murio \cite{Murio08}. We choose this scheme because it is suitable for being transformed into a finite difference  method with non-uniform timesteps and, also, because it allows an easy consideration of discontinuous solutions (see Section \ref{sec:numerical}). For example, this is not true for schemes based in the Gr\"{u}nwald-Letnikov discretization.

The paper is organized as follows. In Section \ref{sec:Alg} the finite difference scheme with non-uniform timesteps is obtained. In Section \ref{sec:Stability} the  stability of the method is analyzed by means of the von-Neumann stability method. Its convergence is proved in Section \ref{sec:convergence}. In Section \ref{sec:numerical}  the method is applied to the problem of evaluating the dispersion of a constant flux of subdiffusive particles stemming from a given point source. The paper ends with some conclusions and remarks in Section \ref{sec:conclusions}.

\section{Fractional difference algorithm with non-uniform timesteps}
\label{sec:Alg}

 As usual in finite difference methods, one starts by considering a mesh in the space-time region where one wants to obtain the numerical estimate $U_j^{(m)}$ of the exact solution $u(x_j,t_m)=u_j^{(m)}$, $(x_j,t_m)$ being the coordinates of the $(j,m)$ node of the mesh. Next, one replaces the continuous operator $\partial$ of the equation $\partial u =F$ one has to solve by a difference operator $\delta$ and a truncation error $R(x,t)$: $\partial u =\delta u +R$.
 For example, in this paper  we replace the continuous operators that define $\partial$ in \eqref{FDEoperator} by
 \begin{align}
 \frac{\partial^\gamma}{\partial t^\gamma} u(x,t) &= \delta_t^\gamma u(x,t)+R_t(x), \\
  \frac{\partial^2}{\partial x^2} u(x,t)&= \delta_x^2 u(x,t)+ R_x(t),
\end{align}
where $\delta_t^\gamma$ and $\delta_x^2$ are the corresponding difference operators. This way one gets
 \begin{equation}
 \left[ \delta_t^\gamma  -  K \delta_x^2\right] u(x,t)=F(x,t)+ R(x,t)
\end{equation}
where $R(x,t)=K R_x(t)-R_t(x)$.  In this case  $ \delta=\delta_t^\gamma -K \delta_x^2$.  Neglecting the truncation term, one gets a difference equation $\delta U=F $ whose solution leads to the finite difference estimate of the exact solution $u(x,t)$ at the mesh points. For a given operator $\partial$ one can consider many different difference operators $\delta$, and hence many different finite difference methods for solving the finite difference equation $\partial u =F$.

We here assume that the spatial size of the mesh  $x_{j+1}-x_j=\Delta x$ is constant, and discretize the Laplacian operator $\partial^2/\partial x^2$ by means of the three-point centred formula so that
\be
\label{Lapladiscre2}
\delta_x^2 u(x_j,t)=\frac{u(x_{j+1},t)-2u(x_{j},t)+u(x_{j+1},t)}{(\Delta x)^2}
\ee
with $R_t(x_j)=O(\Delta x)^2$.
For the fractional Caputo derivative we choose  a discretized operator $\delta_t^\gamma $ that is a generalization of the L1 formula \cite{Oldham} for non-uniform meshes (see the Appendix):
\begin{align}
\label{dCapdiscre2}
\delta_t^\gamma u(x,t_{n})=\frac{1}{\Gamma(2-\gamma)} \sum_{m=0}^{n-1} T_{m,n}^{(\gamma)} \left[ u(x,t_{m+1}^{-})-u(x,t_{m}^{+})\right],
\end{align}
 $R_{t_n}(x)$ being of order $t_n^{1-\gamma}  \max_{0\le m\le n-1}\left(t_{m+1}-t_{m} \right)$. Therefore the truncation error is
\be
\label{Rtrun}
R(x_j,t_n)=O(\Delta x)^2+t_n^{1-\gamma} O\left[ \max_{0\le m\le n-1}\left(t_{m+1}-t_{m} \right)\right] ,
\ee
which goes to zero when the spacing of the time-spatial mesh goes to zero. This means that the method is \emph{consistent} \cite{LeVequeBook}.

Therefore, introducing   \eqref{Lapladiscre2} and \eqref{dCapdiscre2} into  $ \delta=\delta_t^\gamma -K \delta_x^2$, the finite difference equation $\delta U=F$ we have to solve  becomes, after multiplying it by $(t_{n} -t_{n-1})^\gamma$,  $\tilde \partial U=\tilde F$, i.e.,
   \begin{equation}
   \label{metImp1}
   \sum_{m=0}^{n-1} \tilde T_{m,n}^{(\gamma)} \left[ U^{m+1}_{j} -U^{m}_{j}\right]
     -S_n  [U^n_{j+1}- 2 U^n_{j} +U^n_{j-1}]= \tilde F(x_j,t_n)
\end{equation}
where
\begin{align}
S_n&= \Gamma(2-\gamma)  K \frac{ \left(t_{n} -t_{n-1}\right)^\gamma}{(\Delta x)^2},\\
\label{tilTdef}
\tilde T_{m,n}^{(\gamma)}&=  (t_{n}-t_{n-1})^{\gamma}T_{m,n}^{(\gamma)},\\
T_{0,1}^{(\gamma)}& = (t_{1}-t_{0})^{-\gamma},\\
T_{m,n}^{(\gamma)} & =\frac{(t_{n}-t_{m})^{1-\gamma} - (t_{n}-t_{m+1})^{1-\gamma}}{t_{m+1}-t_{m}}, \quad m\le n-1,
\label{Tmn}\\
\tilde F(x_j,t_n)&= (t_{n}-t_{n-1})^{\gamma} F(x_j,t_n).
\end{align}
From here on we use the convention that a tilde over a symbol stands for that symbol multiplied by $(t_n - t_{n-1})^\gamma$.
Reordering  (\ref{metImp1}) and taking into account that $\tilde T_{n-1,n}^{(\gamma)}=1$ [see \eqref{Tmn}] we get the finite difference scheme we were looking for:
\be
\label{ImpliAlg1}
  -S_n \, U^n_{j+1} +(1+2S_n) U^n_{j}-S_n \, U^n_{j-1}=\mathcal{M} U_j^{(n)}+\tilde F(x_j,t_n)
\ee
  $\mathcal{M}$ being the difference operator defined by
\begin{align}
\label{MU}
\mathcal{M} U_j^{(n)} &\equiv U^{n-1}_{j} -   \sum_{m=0}^{n-2} \tilde T_{m,n}^{(\gamma)} \left[ U^{m+1}_{j} -U^{m}_{j}\right] .
\end{align}
Here and in the rest of the paper we adopt the convention that the summation is zero when the lower bound  is larger than the upper bound.

Equation \eqref{ImpliAlg1} can be written in vectorial form  as
\be
\label{ImpliAlg2}
  A U^{(n)}=\mathcal{M} U^{(n)}+\tilde F^{(n)}.
\ee
Comparing this equation with $\tilde \delta U=\tilde F$ one sees that $\tilde \delta =A-\mathcal{M}$.
The present method is implicit because   one has to solve the tridiagonal system (\ref{ImpliAlg1}) in order to get the numerical solution $U^n$.  Fortunately, the system can be efficiently solved by means of the Thomas algorithm  since  $A$  is  a strictly diagonally dominant matrix. For the case of constant timesteps, this method reduces to the one considered in \cite{LiuANZIAM06}.

The operator $\mathcal{M}$ is a kind of difference operator with memory (which comes from the memory, the non-local character, of the fractional derivative) in the sense that its effect on $U$ at time $t_n$, $U^{(n)}$, depends on \emph{all} the previous values $\{ U^{(0)},U^{(1)},\cdots U^{(n-1)} \} \equiv U^{\{n-1\}}$. Where convenient, we will emphasize this fact by writing $\mathcal{M} U^{(n)}$ as $\mathcal{M} \left[ U^{\{n-1\}}\right]$ as this last form has the virtue of making clear that the value of $\mathcal{M} U$ at time $t_n$ is obtained from the values of $U$  at all the previous times $t_{n-1},\cdots,t_0$.  Then, the solution of \eqref{ImpliAlg2} can be written as
\be
\label{eSol}
  U^{(n)}=A^{-1} \mathcal{M} \left[ U^{\{n-1\}}\right]+A^{-1}\tilde F^{(n)}.
\ee

It is interesting to note that, if one discretizes the Caputo derivative at time $t_{n+1}$ and the Laplacian at time $t_n$, then  one straightforwardly gets an \emph{explicit} finite difference scheme
 \begin{align}
 \sum_{m=0}^{n} \tilde T_{m,n+1}^{(\gamma)} \left[ U^{m+1}_{j} -U^{m}_{j}\right]  -S_n  [U^n_{j+1}- 2 U^n_{j} +U^n_{j-1}]=\tilde F(x_j,t_n).
\end{align}
We will not explore this explicit method here because it is unstable if, for a given set of parameters  $\gamma$, $\Delta x$, and $K$, the timesteps are not small enough. This indeed greatly reduces any  advantage of a method with variable timesteps. Nonetheless,  to examine how the non-uniformity of the timesteps affects the regions of stability in the space of parameters would be, from a theoretical perspective, an interesting topic in itself.  For the case of constant timesteps, this explicit method becomes the one considered by the authors in \cite{QuintanaYusteJCND11}.

\section{Stability}
\label{sec:Stability}

In order to check whether the finite difference scheme $\tilde \partial U=\tilde F$ of equations \eqref{metImp1} or \eqref{ImpliAlg1}  is stable, one studies how a perturbed solution $\hat U$ evolves with respect to the reference solution $U$, or, in other words, how (the size of) the perturbation $v=\hat U-U$ evolves in time.
Usually, the initial perturbation $v^{(0)}$ is considered as the difference between the exact initial condition and its computer finite precision representation (round-off error).  Because $\tilde \delta $ is a linear operator, one sees that the perturbation satisfies $\tilde \delta v=0$, i.e., the same equation as $U$ in  \eqref{metImp1}, \eqref{ImpliAlg1}, or \eqref{ImpliAlg2}  when $F=0$:
\be
\label{AvMv}
A v^{(n)}=\mathcal{M} v^{(n)}.
\ee

In order to analyze the stability of the present fractional diffusion difference algorithm, we use the von Neumann-Fourier technique introduced for this kind of schemes in Refs. \cite{YusteSIAM05,YusJCP06}: first one  assumes that $v^{(n)}_{j}$ is described by a discrete Fourier series
\be
\label{vnjFou}
v^{(n)}_{j}=\sum_{q} \xi_q^{(n)} e^{iqj\Delta x},
\ee
where the summation is carried over all the wave numbers $q$ supported by the lattice, and then one analyzes the stability of the complete solution by analyzing the stability of a generic Fourier $q$-mode, say   $\xi_q^{(n)}  e^{iqj\Delta x}$. The rationale for this procedure is that if any mode is stable  then  the complete solution $v^{(n)}_{j}$ written as a superposition of modes is stable too.

Inserting the expression for the generic mode $\xi_q^{(n)}  e^{iqj\Delta x}$ into  \eqref{AvMv} and using the definitions of $A$ and $\mathcal{M} U$ [see \eqref{MU}] one gets
\be
(1+\hat S_n)\xi_q^{(n)}=   \sum_{m=0}^{n-1}\left( \tilde T_{m,n}^{(\gamma)} - \tilde T_{m-1,n}^{(\gamma)}\right)\xi_q^{(m)}
\label{SnTrecur}
\ee
where
\begin{align}
 \hat S_n = 4\sin^2\left(\frac{q\Delta x}{2}\right) S_n .
\end{align}
Here we have taken into account that $\tilde T_{n-1,n}^{(\gamma)}=1$ and  have defined  $\tilde T_{-1,n}^{(\gamma)}\equiv 0$.

Now we depart from our previous \cite{YusteSIAM05,YusJCP06} quick method of analysis in terms of the so-called amplification factor (as it leads to sums involving $T_{m,n}^{(\gamma)}$ we have not been able to evaluate), and we use a procedure similar to the one followed by Murio in Ref. \cite{Murio08}. However, it should be noted that our demonstration is more general (and difficult) as we do not assume, as in Ref. \cite{Murio08},  that the (in general) complex quantity $\xi_q^{(n)}$ is real and positive.

From \eqref{SnTrecur} and taking into account that $\tilde T_{m,n}^{(\gamma)} -\tilde T_{m-1,n}^{(\gamma)}>0$ for any temporal mesh [see Appendix, equation \eqref{Tmngt}], one finds that
\begin{align}
(1+\hat S_n) \left|\xi_q^{(n)}\right|& \le  \sum_{m=0}^{n-1}\left( \tilde T_{m,n}^{(\gamma)} - \tilde T_{m-1,n}^{(\gamma)}\right)\left|\xi_q^{(m)}\right| \le \left| \xi_q^{\{n-1\}}\right|_\text{max} \sum_{m=0}^{n-1}\left( \tilde T_{m,n}^{(\gamma)} - \tilde T_{m-1,n}^{(\gamma)}\right)
\label{SnEst}
\end{align}
where we have defined
\be
\left| \xi_q^{\{n-1\}}\right|_\text{max}\equiv\text{max}\left\{\left| \xi_q^{(0)}\right|,\left| \xi_q^{(1)}\right|,\cdots,\left| \xi_q^{(n-1)}\right| \right\} .
\ee
But   $   \sum_{m=0}^{n-1}\left( \tilde T_{m,n}^{(\gamma)} -\tilde T_{m-1,n}^{(\gamma)}\right)  =1$ [recall that  $\tilde T_{n-1,n}^{(\gamma)}=1$ and  $\tilde T_{-1,n}^{(\gamma)}\equiv 0$] so that
\be
\label{xiEstmax}
\left| \xi_q^{(n)}\right|\leq \left| \xi_q^{\{n-1\}}\right|_\text{max}
\ee
because $(1+\hat S_n)>1$.  Obviously, \eqref{xiEstmax} implies, a fortiori, that
\be
\label{xiEst}
\left| \xi_q^{(n)}\right|\leq \left| \xi_q^{(0)}\right|
\ee
for all  $n$.
Therefore, the perturbation of the generic mode remains smaller than or equal to its initial perturbation $\left| \xi_q^{(n)}\right|\leq \left| \xi_q^{(0)}\right|$. This means that the present difference method is unconditionally stable.

Often the stability condition is expressed  in terms of the 2-norm (or Euclidean norm) of the perturbation:
$\left\|v^{(n)}\right\|_{2}^2\equiv \Delta x\, \sum_{j} \left|  v_j^{(n)} \right|^2$. One can transform the stability \eqref{xiEst} to this form by means of  Parseval's relation, $\left\|v^{(n)}\right\|_{2}^2= \sum_{q} \left|\xi_q^{(n)}\right|^2 $ \cite{MortonMayers}. Using this expression \eqref{xiEst} becomes
\begin{align}
\label{vnlvn1}
\left\|v^{(n)}\right\|_{2}  \le  \left\|v^{(0)}\right\|_{2},
\end{align}
which implies that the norm of the perturbation remains bounded by the initial perturbation:
$ \left\|v^{(n)}\right\|_{2}   \le  \left\|v^{(0)}\right\|_{2} $.

The solution of \eqref{AvMv} can be written  as
\be
\label{vSol}
  v^{(n)}=A^{-1}\mathcal{M} \left[ v^{\{n-1\}}\right]
\ee
so that the stability condition implied by \eqref{vnlvn1} means that the operator $A^{-1}\mathcal{M}$ has the key property
\begin{align}
\label{Sesta}
\left\|A^{-1}\mathcal{M}\left[ v^{\{n-1\}}\right] \right\|_{2}  \le  \left\|v^{(0)}\right\|_{2}.
\end{align}

\section{Consistency and stability implies convergence}
\label{sec:convergence}

The basic property a difference scheme $\delta U=F$ should have is that its approximate solutions should converge towards the exact solution of $\partial U=F$ when the size of the spatiotemporal discretization goes to zero. In this case it is said that the method is convergent \cite{LeVequeBook,MortonMayers}. In this section we show that for a whole class of fractional difference algorithms, similarly to the case for the usual difference algorithms, the consistency and stability of the difference scheme implies its convergence.  (This result can be seen as the extension  of the sufficient condition of  Lax's equivalence theorem  to fractional equations of the form $\partial u=F$, or as the fractional counterpart of what is sometimes called fundamental theorem of finite difference methods \cite{LeVequeBook}.)

 Let us define $e_j^{(k)}$ as the difference between the numerical and exact solution at the point $(x_j,t_m)$:  $e_j^{(k)}= U_j^{(k)}-u_j^{(k)}$. By definition $\partial u=\delta u+R$,  $ \delta U=0$ and $\partial u=0$, so that, because these operators are \emph{linear},  $\tilde \delta e= \tilde R$.  But  $\tilde \delta=A-\mathcal{M}$ [see following  equation \eqref{ImpliAlg2}], and then
\be
A e^{(n)} = \mathcal{M} \left[ e^{\{n-1\}}\right]+   \tilde R^{(n)} \Leftrightarrow e^{(n)} =A^{-1}\mathcal{M}\left[ e^{\{n-1\}}\right]+ A^{-1} \tilde R^{(n)}.
\ee
Therefore
\be
\label{en4}
\left\|e^{(n)}\right\|_{2} \le \left\| A^{-1}\mathcal{M}\left[ e^{\{n-1\}}\right] \right\|_{2}+ \left\|  A^{-1} \right\|_{2} \, \left\| \tilde  R^{(n)}\right\|_{2}
\ee
and,  because the  method is \emph{stable},   \eqref{Sesta} holds  and \eqref{en4} becomes
\be
\label{enfin}
\left\|e^{(n)}\right\|_{2} \le \left\|  e^{(0)} \right\|_{2}+ \left\|  A^{-1} \right\|_{2} \, \left\| \tilde  R^{(n)}\right\|_{2}.
\ee
If, when the size of the discretization mesh goes to zero, (i) $\left\|\bar R^{(n)}\right\|_{2}$ goes to zero (i.e., if the method is \emph{consistent}) and (ii)  $\left\|  A^{-1} \right\|_{2}$ is bounded, then \eqref{enfin} implies that the error   $\left\|e^{(n)}\right\|_{2}$ goes to zero. i.e., the method is convergent.
 In particular, the finite difference method of this paper is convergent because, (i) the method is consistent [see following equation \eqref{Rtrun}] and (ii) $\left\|  A^{-1} \right\|_{2}$ is always bounded since  $A$  is a real symmetric tridiagonal Toeplitz matrix with all its eigenvalues greater than 1 \cite{LeVequeBook}.
Note that the present demonstration is quite general as only generic properties of the operators $\mathcal{M}$ and $A$ appearing in the finite difference equation (linearity, boundedness, stability) are required. In fact this demonstration is \emph{parallel} to the standard procedure to prove this result for non-fractional diffusion problems \cite{AllaireBook}.  This way, one realizes that some aspects of fractional difference methods can be included into the realm of standard non-fractional finite difference theory. In fact, fractional difference methods could be seen as particular cases of multilevel schemes where the number of levels \emph{increases} with time.

\section{Dispersion of a flux of subdiffusive particles: Numerical results}
\label{sec:numerical}

In this section we test the finite difference scheme  (\ref{metImp1}) by solving  a problem that describes the dispersion of a constant flux of anomalous subdiffusive particles appearing at  a given place. In particular,  in  our problem the flux is 1 as one particle is released every unit time at times $t=0,1,2\ldots$ from a source of particles placed at $x=0$ in a  one-dimensional infinite medium. We shall assume that the probability density $u(x,t)$  of finding the particle at position $x$ at time $t$ follows the fractional diffusion equation $\partial u=0$. The release of a particle at time $n$ is described by the introduction of  the  Dirac  delta function $\delta(x)$ at that time [see following equation \eqref{CaputoRL2b}].
The probability distribution associated with each particle is just the propagator (or Green's function) $G(x,t-k)$ of the problem, i.e., the solution of $\partial u=0$ in the unbounded space when there is only a single Dirac delta function that appears at $t=k$. This propagator can be explicitly written in terms of  Fox's $H$ function  \cite{MetzlerKlafterPhysRep,MathaiSaxena}:
\begin{equation}
\label{propag} G(x,t)=\frac{1}{\sqrt{4\pi K_\gamma t^\gamma}}\,
H^{10}_{11}\left[\frac{|x|}{\sqrt{ K  t^\gamma}}
\left|\begin{array}{l}{(1-\gamma/2,\gamma/2)}\\[1ex]{(0,1)}\end{array}\right.\!\!\right].
\end{equation}
Finally, given the linear nature of the problem, its exact solution is just a superposition of Green functions:
\be
u(x,t)=\sum_{k=0}^{\lfloor t\rfloor} G(x,t-k)
\ee
where ${\lfloor t\rfloor}$ is the floor function.  In the formalism of \eqref{FDEdef},  the function $F$ of this problem is
 \be
F(x,t)=\delta(x) \sum_{k=0}^{\lfloor t\rfloor}  \delta^{(\gamma)}(t-k)
\ee
 with
   \be
  \delta^{(\gamma)}(t-k)   \equiv \frac{\partial^\gamma}{\partial t^\gamma} \Theta(t-k),
\ee
 $\Theta$ being the Heaviside step function.

 The problem just presented, with the added feature that the particles may react and disappear, describes the development of morphogen  gradients (a key concept in developmental biology) when the anomalous diffusion of the morphogens is accounted for by means of a CTRW model \cite{YusAbaLinPRE11}. Because knowledge of the long-time distribution of morphogens in the medium is especially important, numerical methods that can use large timesteps (to go fast in time) and small timesteps (to take account of regions of large variability when the particles are released) would seem to be particularly convenient.

Regarding the boundary conditions, although it is of course not possible to use boundary conditions at infinity in numerical calculations,  one can solve the problem by means of some convenient boundary conditions, say absorbing boundary conditions  $u(x=-L/2,t)=u(x=L/2,t)=0$, so far away from the source at $x=0$ that,  for a finite $t$, their impact on the solution is negligible in the region of interest which is, by construction, far away from the boundaries.  In our computations, where  we have chosen $K=1$ and $\gamma=1/2$ and the largest value of time is $t=2$,  we used $L=20$. For this choice the exact value of the solution at the boundaries is  $u(\pm 10,2)\simeq 5\times 10^{-5}$, which can be safely approximated by 0  in our numerical calculations (as we will show shortly, the errors are well above  $10^{-4}$).
The nodes of the spatial mesh we used were placed at  $x_j=j\,\Delta x$ with $j=-N/2, -N/2+1,\cdots, N/2$,     $N=100$,  and $\Delta x=0.2$. Finally, every Dirac  delta function  $\delta(x)$ was approximated by its discretized version:
 \begin{equation}
\delta(x)\simeq \left\{\begin{array}{ll}
 1/\Delta x , & \quad j=0,\\
 0 ,      &  \quad   j\neq 0.
 \end{array}
 \right.
 \label{deltaApprox}
 \end{equation}
Note that the  discontinuity of the solution due to the periodic appearance of Dirac  delta functions at times $t=k$ ($k$ integer) has yet to be worked out. Although in section \ref{sec:Alg} we only consider the case in which the solution is continuous, discontinuities can be easily handled by our algorithm taking into account  \eqref{dCapdiscre2}, which is valid even for discontinuous functions at the discretization times $t_m$.
This way, keeping the symbol $U_j^{(m)}$ for the estimate of $u(x_j,t_m^-)$ and denoting by $V_j^{(m)}$   the estimate of $u(x_j,t_m^+)$,  \eqref{metImp1}  becomes
   \begin{equation}
   \label{metImp9}
   \sum_{m=0}^{n-1} \tilde T_{m,n}^{(\gamma)} \left[ U^{m+1}_{j} -V^{m}_{j}\right]
     -S_n  [U^n_{j+1}- 2 U^n_{j} +U^n_{j-1}]= \tilde F(x_j,t_n) .
\end{equation}
The algorithm in this case  can be written as in  \eqref{ImpliAlg2}  but now with
 \begin{align}
\mathcal{M} U_j^{(n)} &\equiv V^{n-1}_{j} - \sum_{m=0}^{n-2} \tilde T_{m,n}^{(\gamma)} \left[ U^{m+1}_{j} -V^{m}_{j}\right] .
\end{align}
 In particular, the appearance of a Dirac  delta function at $x=0$  at time $t=k$ is handled by simply making $V_0^{(n)}=1/\Delta x+U_0^{(n)}$ whenever $t_n=k$.

\begin{figure}
\begin{center}
\includegraphics[width=0.7\columnwidth]{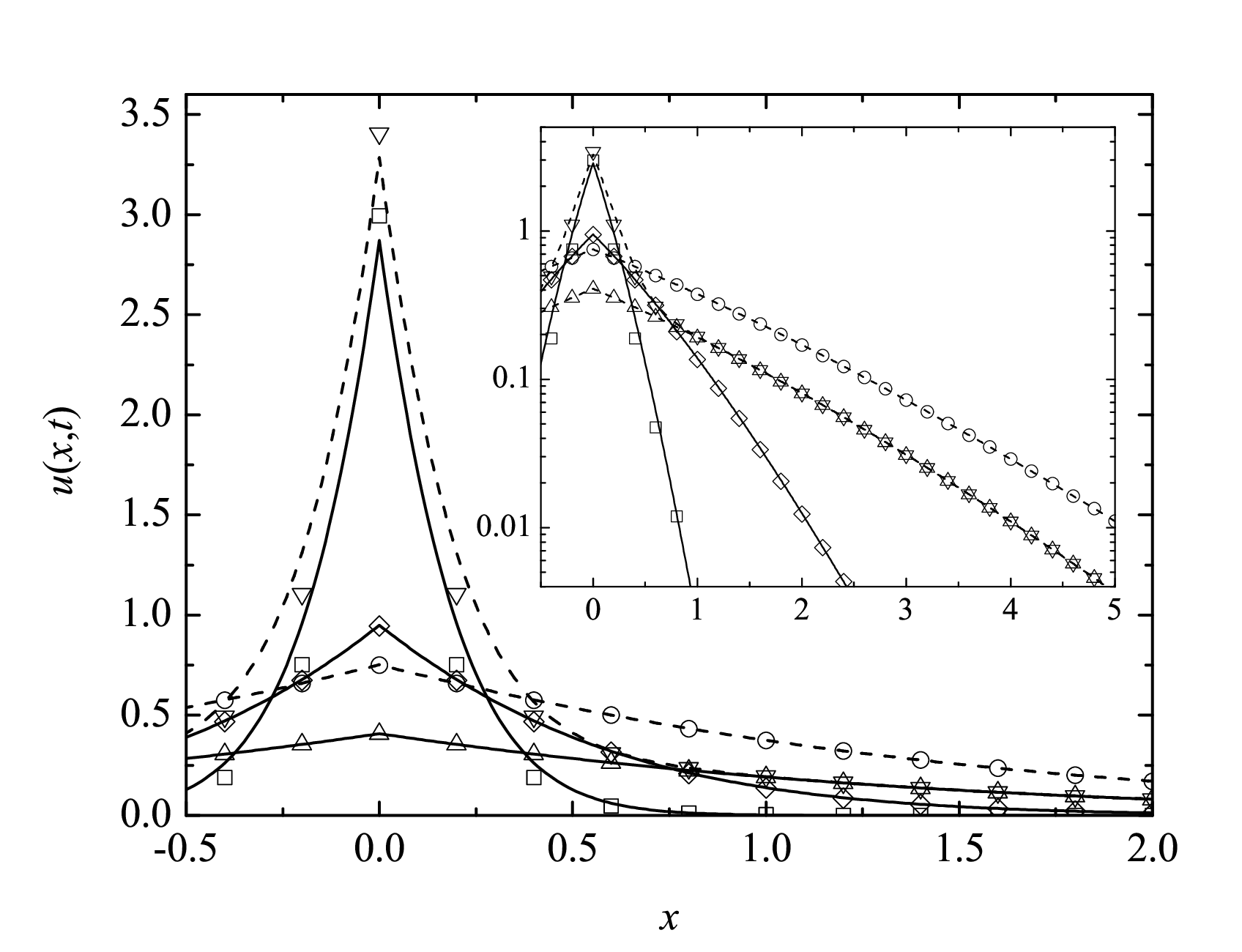}
\caption{ Exact solution $u(x,t)$ (lines) and numerical solutions $U_j^{(n)}$ with variable timesteps given by (\ref{dtn})  (symbols) for the problem described in the main text with (from top to bottom at $x=0$) $t=1.0004$ (66 timesteps, down triangles),   $t=4.08\times 10^{-4}$ (one timestep, squares), $t=0.034$ (10 timesteps, diamonds),  $t=2.0$ (141 timesteps, circles), and  $t=1$ (65 timesteps, up triangles). The lines for $t>1$ are dashed. Inset: logarithmic scale and tails of the solutions.}
\label{soluxt}
\end{center}
\end{figure}

\begin{figure}
\begin{center}
\includegraphics[width=0.7\columnwidth]{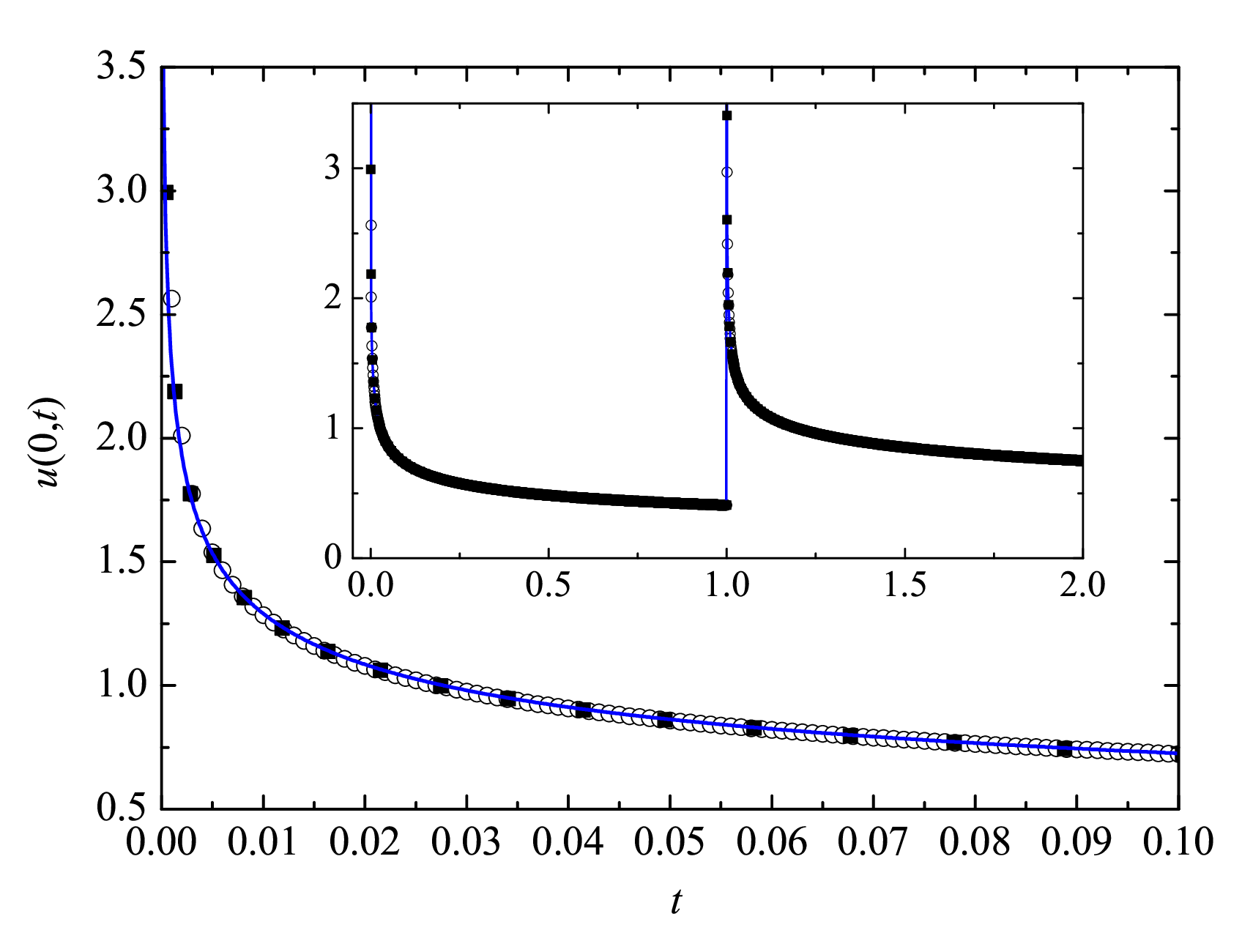}
\caption{ Exact solution $u(0,t)$ (line) and numerical solution $U_0^{(n)}$ (symbols)  at the origin $x=0$ for the problem described in the main text for fixed timesteps with $t_{m+1}-t_m=0.001$ (circles) and for variable timesteps given by (\ref{dtn}) (solid squares) in the time interval $0\le t\le 0.1$. Inset: solutions in the whole time interval. In this panel only 1 of every 20 points for the case with fixed timesteps is plotted.}
\label{solu}
\end{center}
\end{figure}

In Fig. \ref{soluxt} we compare  exact solutions and  numerical solutions obtained by means of the present implicit method when the non-uniform discretization is given by this adaptive timestep:
\be
\label{dtn}
t_{m+1}-t_m= \text{min}\left[10^{-4}\coth\left[|g(t_m)|/1000\right], 0.02\right].
\ee
Here $g(t_m)=(U_{-1}^{m}-2U_0^{m}+U_{1}^{m})/(\Delta x)^2$ is the three-point centred estimate of the spatial second derivative  at $x=0$ and time $t_m$: $g(t_m)\simeq \partial^2 u(x,t_m)/\partial x^2|_{x=0}$, i.e., an estimate of the curvature of the solution at the place where it changes most  abruptly. The agreement is excellent. In fact, even at $x=0$, the place where the numerical solution is worst, the results are quite good (see Fig.~\ref{solu}). The relatively large error for short times after the appearance of the Dirac  delta function is due to its crude approximation by (\ref{deltaApprox}).

\begin{figure}
\begin{center}
\includegraphics[width=0.7\columnwidth]{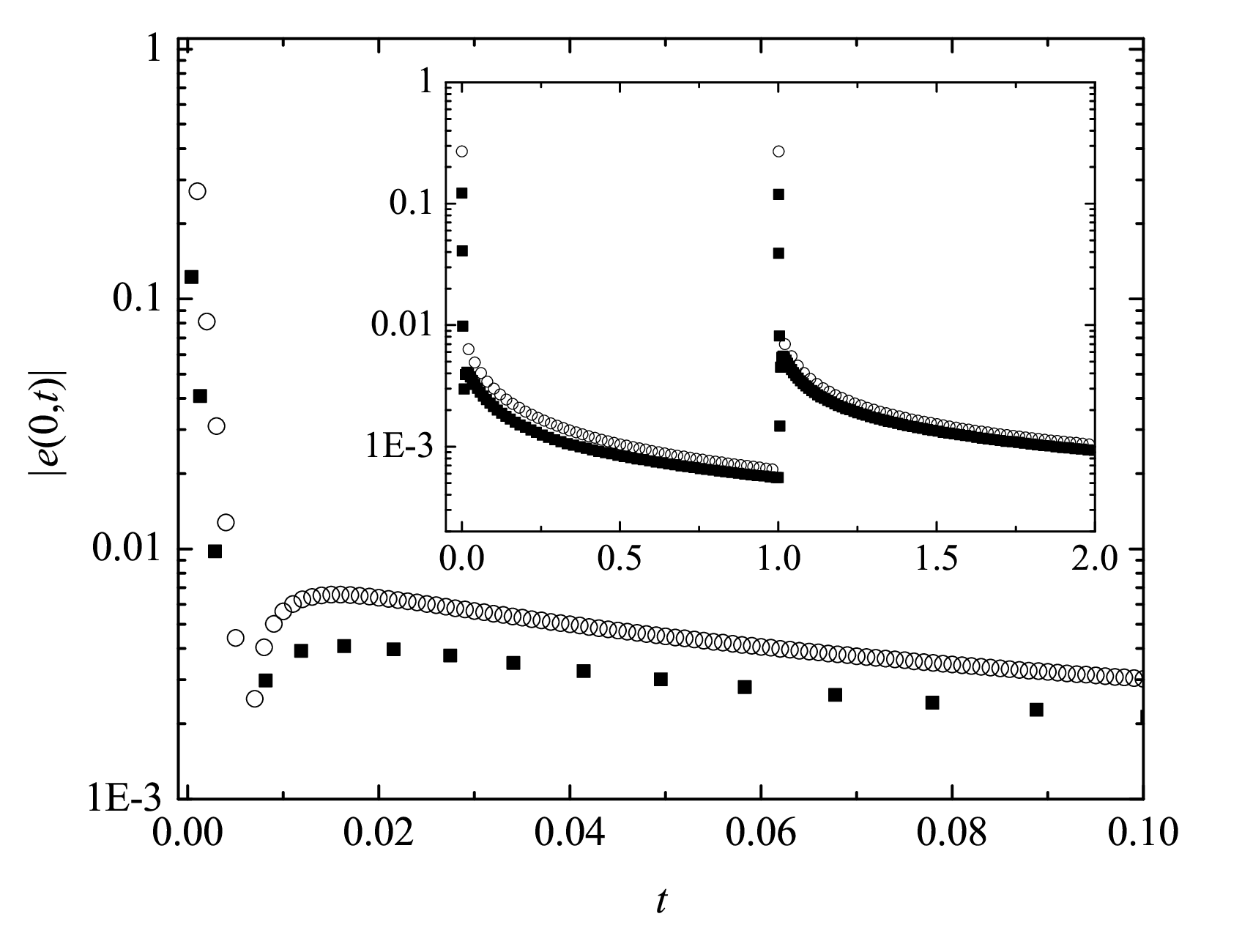}
\caption{\label{error}  Absolute difference between the  exact solution  and the numerical solution at $x=0$, $|e(0,t)|$, for the problem described in the main text for fixed timesteps with $t_{m+1}-t_m =0.001$ (circles) and variable timesteps given by (\ref{dtn}) (squares) in the time interval $0\le t\le 0.1$. Inset: error in the whole time interval. In this small panel only 1 of every 20 points for the case with fixed timesteps is plotted.
}
\end{center}
\end{figure}

Although formula~(\ref{dtn}) for the temporal spacing $t_{m+1}-t_m$ is indeed arbitrary, it has some convenient features: it has prefixed minimum and maximum values (0.0001 and 0.02, respectively) and $t_{m+1}-t_m$  increases (decreases) when the curvature at $x=0$, $ \partial^2 u(x,t_m)/\partial x^2|_{x=0}$ decreases (increases).  This guarantees a more thorough  description of the solution where necessary, namely, those regions where the solution changes most quickly.  For example, when (\ref{dtn}) is used, 5 timesteps are employed in the tiny temporal region between $t=0$  and $t=0.01$ (a region where the solution changes abruptly; see Fig.~\ref{solu}); however, only 59 more timesteps are necessary to cover the comparatively much larger region between $t=0.01$ and $t=1$.  Note that this makes the procedure with variable timesteps as least as accurate as the  method with fixed timesteps  $t_{m+1}-t_m=0.001$ in the difficult region near $t=0$ (see Fig.~\ref{error}). However this is not reached at at the expense of computational cost: compare the $64$ steps [17 steps] required by the variable timesteps method to reach $t=1$ [$t=0.1$] with the $1000$ steps [100 steps] required when $t_{m+1}-t_m=0.001$. This difference is far from minor because, from  (\ref{metImp1}), one sees that the computation time required to calculate the solution after $n$ timesteps increases roughly as $n^2$. In fact, in our actual calculations we found that the CPU time spent to evaluate the solution at $t=1$  with $t_{m+1}-t_m=0.001$ was well above one hundred times! longer than when (\ref{dtn}) was used.

Finally, it should be noted that for a standard (non-fractional) diffusion equation a time discretization based on the estimate $g(t)$ of $\partial^2 u/\partial x^2$  (which is the procedure employed in this paper) is equivalent to a time discretization based on the estimate of the (local) derivative  $\partial  u/\partial t$.  This is due to the fact that both quantities are related by the diffusion equation  so that to estimate $\partial^2 u/\partial x^2$ is tantamount to estimate $\partial  u/\partial t$.
Similarly, for the fractional diffusion equation  \eqref{FDEdef},  a time discretization based on the estimate $g(t)$ of $\partial^2 u/\partial x^2$ (Laplacian)  is equivalent to a discretization based on the estimate of the (fractional) derivative $\partial  u^\gamma/\partial t^\gamma$. However, because the Laplacian is a local operator (whereas fractional derivatives are non local operators) the easiest way to estimate the size of $\partial  u^\gamma/\partial t^\gamma$ at every timestep is by calculating $g(t)$.  Note, finally, that for this kind of fractional problem the local time derivative is not fully appropriate to guess the behavior of the solution (and, accordingly, the size of the timesteps) because the behavior of the solution depends on its presents and \emph{past} values due to the non-local nature of the fractional time derivative.  In summary, the rationale behind our choice of   $g(t)$  as the parameter to gauge the size of the timesteps is the consideration that the temporal change of the solution depends on the size of the fractional time derivative which, in turns, can be estimated from the value of the Laplacian by using the fractional diffusion equation \eqref{FDEdef}.

\section{Conclusions and final remarks}
\label{sec:conclusions}

An implicit finite difference method with non-uniform timesteps for fractional diffusion equations has been described. The method allows the temporal mesh to be adapted to the behaviour of the solution by sampling thoroughly those regions where the variation of the solution is large --so as to maintain the accuracy of the method-- but sparsely wherever the solution varies mildly --and so advancing fast in time and reducing the number of timesteps needed in the computation. This last feature is especially convenient when numerically solving  fractional problems since the computational cost in this case scales as the \emph{square} of the number of timesteps. Compare this with the typical linear dependence for non-fractional diffusion problems.

The finite difference method discussed in this paper was obtained from the fractional diffusion equation in the Caputo form by discretizing (i) the Laplacian by means of the standard three-point centred formula and (ii) the Caputo derivative with a generalization of the L1 formula to non-uniform meshes.
The corresponding truncation error is of second order in the size of the spatial mesh and of first order in the size of the temporal mesh.  A von-Neumann stability analysis was employed to prove that the method is unconditionally stable. Also it was shown in a quite general form that the consistency and stability of the method implies its convergence.  In fact, this implication, sometimes called the fundamental theorem of finite difference methods,  it is not restricted to the particular finite difference scheme considered in this paper, but it is a quite general result that depends only on some general properties of the fractional finite difference scheme.

The present method can be generalized and extended in several ways. First, one could employ other different non-uniform discretizations of the fractional derivative, say of higher-order accuracy, or even consider non-uniform discretization in the space.
Also one could explore criteria or formulas for efficiently choosing the size of the timesteps (adaptive methods). Finally, another interesting avenue to explore is the application of this kind of method to other classes of fractional equation, such as fractional diffusion-wave equations, fractional Fokker-Planck equations, and fractional equations with Riemann-Liouville derivatives.

\section*{Acknowledgments} \noindent
 Financial support from the Ministerio de Educaci\'on  y Ciencia through Grant No. FIS2010-16587 (partially financed with FEDER funds) and by the Junta de Extremadura (Spain) through Grant No. GRU10158 is gratefully acknowledged.

\appendix
\section{Non-uniform L1 formula for the Caputo derivative}

We start by rewriting the Caputo derivative defined in   \eqref{CaputoRL2b}:
\begin{align}
\left. \frac{d^\gamma y(t)}{d t^\gamma}\right|_{t=t_{n}}
&= \frac{1}{\Gamma(1-\gamma)} \sum_{m=0}^{n-1}
 \int_{t_{m}}^{t_{m+1}} d\tau  \frac{1}{(t-\tau)^{\gamma}} \frac{d y(\tau)}{d\tau}
\end{align}
with $t_0=0$ and $n\ge 1$. Proceeding as is standard in the so-called L1 method \cite{Oldham} the ordinary derivative inside the integral with lower limit $t_m$ is approximated by
 \be
 \frac{d y(t)}{d t} = \frac{y(t_{m+1}^-)-y(t_{m}^+)}{t_{m+1}-t_{m}} + O(t_{m+1}-t_{m}) ,  \quad t_{m}\le t\le t_{m+1}
 \ee
so that
 \begin{align}
 \label{dytg1}
\left. \frac{d^\gamma y(t)}{d t^\gamma}\right|_{t=t_{n}}
&= \frac{1}{\Gamma(1-\gamma)} \sum_{m=0}^{n-1} \frac{y(t_{m+1}^-)-y(t_{m}^+)}{t_{m+1}-t_{m}}
 \int_{t_{m}}^{t_{m+1}} d\tau   (t-\tau)^{-\gamma} +R(t_n)
\end{align}
where $y(t_{m}^-)$ and  $y(t_{m}^+)$ are the values of  $y(t)$ when $t$ goes to $t_m$ from the left and right, respectively, and
\be
\label{Req}
 R(t_n)=\frac{1}{\Gamma(1-\gamma)}  \sum_{m=0}^{n-1} R_m   \int_{t_{m}}^{t_{m+1}} d\tau   (t-\tau)^{-\gamma}
\ee
is the temporal truncation error.
The evaluation of the integral of \eqref{dytg1} is straightforward and this equation becomes
\begin{align}
\label{dCapdiscre1}
\left. \frac{d^\gamma y(t)}{d t^\gamma}\right|_{t=t_{n}}
&=
\frac{1}{\Gamma(2-\gamma)} \sum_{m=0}^{n-1} T_{m,n}^{(\gamma)} \left[ y(t_{m+1}^{-})-y(t_{m}^{+})\right] + R(t_n)
\end{align}
where $T_{m,n}^{(\gamma)}$ is given is \eqref{Tmn}. Therefore, from \eqref{Req} one sees that
\be
 \left| R(t_n) \right| \le \frac{1}{\Gamma(1-\gamma)}   \max_{0\le m\le n-1} \left|R_m\right|   \int_{0}^{t_{n}} d\tau   (t-\tau)^{-\gamma}
\ee
so that $R(t_n)\le C \max_{0\le m\le n-1} \left|R_m\right|$,  $C$ being a constant of order $t_n^{1-\gamma}$, i.e.,
\be
R(t_n) \le  t_n^{1-\gamma} O\left[\max_{0\le m\le n-1}\left(t_{m+1}-t_{m} \right)\right].
\ee

Finally, we shall prove that
\be
\label{Tmngt}
T_{m,n}^{(\gamma)}>T_{m-1,n}^{(\gamma)},
\ee
a result we used in  \eqref{SnEst}. From the definition of $T_{m,n}^{(\gamma)}$ and making the change $z_m=t_n-t_m$, one sees that
\be
T_{m,n}^{(\gamma)}=\frac{z_m^{1-\gamma}-z_{m+1}^{1-\gamma}}{z_{m}-z_{m+1}}.
\ee
Let us define $f(z)=z^{1-\gamma}$. To prove \eqref{Tmngt}  is equivalent to proving that
 \be
 \frac{f(z_{m+1})-f(z_m)}{z_{m+1}-z_{m}} <\frac{f(z_{m})-f(z_{m-1})}{z_{m}-z_{m-1}} ,
\ee
which is true for any monotonous growing function $f(z)$ with monotonous decreasing derivative $f'(z)$ such as $f(z)=z^{1-\gamma}$.


\begin{thebibliography}{99}

\bibitem{Oldham} K.B. Oldham, J. Spanier, The Fractional Calculus ( Academic Press, New York, 1974).

\bibitem{HilferEd} R. Hilfer (ed.), Applications of Fractional Calculus in Physics (World Scientific, Singapore, 2000).

\bibitem{MetzlerKlafterPhysRep} R. Metzler,  J. Klafter, Phys. Rep. 339.
    (2000) 1.

\bibitem{PhysToday} I. M. Sokolov, J. Klafter, A. Blumen,
   Phys. Today 55 (2002) 48.

\bibitem{MetzlerKlafterJPA04} R.~Metzler, J.~Klafter,  J. Phys. A-Math. Gen. 37   (2004)  R161


\bibitem{Kilbas:06} A.A. Kilbas, H.M. Srivastava, J.J. Trujillo, Theory and Applications of Fractional Differential Equations ( Elsevier, Amsterdam, 2006).

\bibitem{Anotrans:08} R. Klages, G. Radons, I.M. Sokolov (eds.),
 Anomalous Transport: Foundations and Applications
(Elsevier, Amsterdam, 2008).

\bibitem{Magin:08} R.L. Magin, O. Abdullah, D. Baleanu, X.J. Zhou,
 J. Mag. Reson. 190   (2008) 255.




\bibitem{DixVerkman08} J.A. Dix, A.S. Verkman, Ann Rev. Biophys.  37 (2008) 247.


\bibitem{HenryLanglands08} B.I. Henry, T.A.M. Langlands, S. Wearne,  Phys. Rev. Lett. 100   (2008)  128103.


\bibitem{BarMetKlaPRE00} E. Barkai, R. Metzler, J. Klafter,  Phys. Rev. E 61  (2000)
    132.


\bibitem {YusAcePhysica04} S.B. Yuste, L. Acedo,
  Physica A  336   (2004) 334.

\bibitem{Sayed} A.M.A. El-Sayed,  M. Gaber,
 Phys. Lett. A 359   (2006) 175.

\bibitem{Jafari}
 H. Jafari, S. Momani, Phys. Lett. A 370   (2007) 388.



\bibitem{Ray}  S.S. Ray, Phys. Scripta 75   (2007) 53.


\bibitem{GorenfloNLD02} R. Gorenflo,  F. Mainardi, D. Moretti, P. Paradisi, Nonlinear Dynam. 29   (2002) 129.


\bibitem{YusteSIAM05} S.B. Yuste and L. Acedo, SIAM J. Numer. Anal. 42   (2005) 1862.


\bibitem{YusJCP06}  S.B. Yuste, J. Comput. Phys. 216 (1) (2006) 264


 \bibitem{LynchJCP03}  V. E. Lynch, B. A. Carreras, D. del-Castillo-Negrete, K. M. Ferreira-Mejias, H. R. Hicks, J. Comput. Phys. 192   (2003) 406.

\bibitem{MeerschaertTadjeranJCP4} M.M. Meerschaert and C. Tadjeran, J.
    Comput. Appl. Math. 172  (2004) 65.


\bibitem{Chen07} C. M. Chen, F. Liu, I. Turner and V. Anh, J. Comput. Phys.  227   (2007) 886.

\bibitem{Sun08} Z. Z. Sun, X. Wu, Appl. Numer. Math. 56   (2006) 193.



\bibitem{Podlubny09} I. Podlubny, A. V. Chechkin, T. Skovranek, Y. Chen, B. M. Vinagre, J. Comput. Phys. 228   (2009) 3137.

\bibitem{Cui09} M. Cui, J. Comput. Phys. 228   (2009) 7792.

\bibitem{BrunnerJCP10} H. Brunner, L. Ling, M. Yamamoto,
 J. Comput. Phys. 229   (2010)  6613.


\bibitem{SkovranekFDA10} T. Skovranek, V. V. Verbickij, Y. Tarte, I. Podlubny,
 Discretization of fractional-order operators and fractional differential equations on a non-equidistant mesh, Article no. FDA10-062, in: I. Podlubny, B. M.
 Vinagre Jara, YQ. Chen, V. Feliu Batlle, I. Tejado Balsera (Eds.), Proceedings of FDA10, The 4th IFAC Workshop Fractional Differentiation and its Applications,
 Badajoz ( 2010) p. 18.


\bibitem{Mustapha11} K. Mustapha, W. McLean, Numer. Algorithms 56   (2011)     159.


\bibitem{DiethelmND02}
 K. Diethelm, N.J. Ford, A.D. Freed,  Nonlinear Dynam.  29   (2002) 3.

\bibitem{DiethelmCMAME05} K. Diethelm, N.J. Ford, A.D. Freed, Yu. Luchko, Comp. Method. Appl. M.  194   (2005) 743.

\bibitem{FordNA01} N.J. Ford, A.C. Simpson, Numer. Algorithms 26   (2001) 333.

\bibitem{Murio08} D.A. Murio, Comput. Math. Appl. 56   (2008) 1138.

\bibitem{YusteQuintanaPS09} S.B. Yuste, J. Quintana-Murillo,  Phys. Scripta  T136  (2009) 014025.

\bibitem{LiuANZIAM06} F. Liu, P. Zhuang, V. Anh,  I. Turner, ANZIAM J.  47 (2006) C48.

\bibitem{LeVequeBook} R.J. LeVeque, Finite Difference Methods for Ordinary and Partial Differential Equations: Steady-State and Time-Dependent Problems (SIAM,
    Philadelfia, 2007).

\bibitem{QuintanaYusteJCND11}  J. Quintana-Murillo, S.B. Yuste,
 J. Comput. Nonlin. Dyn. 6  (2011)  021014.

\bibitem{MortonMayers} K.W. Morton,  D.F. Mayers, Numerical Solution of Partial Differential Equations,
  (Cambridge University Press, Cambridge, 1994).

\bibitem{AllaireBook}  G. Allaire,  Numerical Analysis and Optimization: An Introduction to Mathematical Modelling and Numerical Simulation (Oxford University
    Press, New York, 2007).

\bibitem{MathaiSaxena} A.M. Mathai, R.K. Saxena, The H-function with Applications in Statistics and Other Disciplines (Wiley, New York, 1978).

\bibitem{YusAbaLinPRE11} S.B. Yuste, E. Abad, K. Lindenberg,  Phys. Rev. E 82   (2010) 061123.


\end{thebibliography}
\end{document}